\documentclass[11pt,reqno]{amsart}

\usepackage{amssymb}
\usepackage{amscd}
\usepackage{amsfonts}
\usepackage{mathrsfs}
\usepackage{setspace}
\usepackage{version}
\usepackage{mathtools}
\usepackage{textcomp}
\usepackage[english]{babel}

\numberwithin{equation}{section} 
\theoremstyle{plain}
\newtheorem{theorem}[subsection]{Theorem}
\newtheorem{lemma}[subsection]{Lemma}
\newtheorem{definition}[subsection]{Definition}
\newtheorem{remark}[subsection]{Remark}
\newtheorem{proposition}[subsection]{Proposition}

\DeclareMathOperator*{\esssup}{ess.sup\,}
\DeclareMathOperator*{\essinf}{ess.inf\,}
\newcommand{\F}{\mathcal{F}}
\newcommand{\EE}{\mathcal{E}}
\renewcommand{\leq}{\leqslant}
\renewcommand{\geq}{\geqslant}

\newcommand\E{\mathbb{E}}
\newcommand\R{\mathbb{R}}
\newcommand\N{\mathbb{N}}
\renewcommand\P{\mathbb{P}}

\newcommand\Q{\mathbb{Q}}

\newcommand{\vertiii}[1]{{\left\vert\kern-0.25ex\left\vert\kern-0.25ex\left\vert #1 
    \right\vert\kern-0.25ex\right\vert\kern-0.25ex\right\vert}}

\begin{document}

\title[]{Duality and stable compactness in Orlicz-type modules}

\author{Jos\'e Orihuela}
\address{Departamento de Matem\'aticas, Universidad de Murcia}
\email{joseori@um.es}

		\author{Jos\'e M.~Zapata}
	\address{Departamento de Estadística e Investigación Operativa, Universidad de Murcia}
	\email{jmzg1@um.es}

\thanks{The authors are grateful to an anonymous referee for a careful review
of the manuscript and valuable comments which have improved the presentation of the results.}


\begin{abstract} 
 Orlicz-type modules are module analogues of classical Orlicz spaces.  
 We study duality and stable compactness in Orlicz-type modules. 
 We characterize the conditional K\"{o}the dual of an Orlicz-type module as the space of all  $\sigma$-order continuous module homomorphisms.   
 We find an order continuity criterion for stable compactness in Orlicz-type modules.  
 As an application, we obtain a robust representation result for conditional risk measures on Orlicz spaces. 
  
\smallskip
\noindent \textit{Orlicz-type module, stable compactness, duality, conditional risk measures}
\end{abstract}

\maketitle

\setcounter{tocdepth}{1}

\section{Introduction}
Motivated by questions in mathematical finance Filipovi\'{c}, Kupper and Vogelpoth in~\cite{kupper03,kupper11,vogelpoth2009phd} developed an analogue of convex analysis for modules over the ordered ring $L^0$ of almost surely equal random variables.   
In particular, they studied topological $L^0$-modules of $L^p$ and Orlicz type. 
In the present paper, we focus on some questions related to Orlicz-type modules. 
We denote by $L^\phi_\F(\EE)$ the Orlicz-type $L^0(\F)$-module  associated with the Young function $\phi\colon[0,+\infty)\to[0,+\infty]$, where $(\Omega,\mathcal{E},\P)$ is a probability space and $\F\subset\EE$ is a sub-$\sigma$-algebra. Formally, we define $L^\phi_\F(\EE) = L^0(\F) L^\phi(\EE)$, where $L^\phi(\EE)$ is the Orlicz space corresponding to the Young function $\phi$. 

Next, we explain the main contributions of the present work.

First, by introducing a module analogue of the K\"{o}the dual of the Orlicz-type module $L^\phi_\F(\EE)$,  we prove that every $\sigma$-order continuous module homomorphism $\mu\colon L^\phi_\F(\EE)\to L^0(\F)$ corresponds to a unique element $y\in L^{\phi^\ast}_\F(\EE)$ such that $\mu(x)=\mathbb{E}[xy|\F]$ for all $x\in L^{\phi^\ast}_\F(\EE)$, where $\phi^\ast$ is the conjugate Young function of $\phi$.

Second, the duality of modules, defined by  $(x, y) \mapsto \mathbb{E}[xy|\F]$, induces the stable weak topology $\sigma_s(L^{\phi^\ast}_{\F}(\EE),L^{\phi}_{\F}(\EE))$. This topology is an instance of the locally $L^0(\F)$-convex topologies introduced in~\cite{kupper03}. We emphasize that locally $L^0(\F)$-convex topologies, as proven in~\cite{jamzap}, are anti-compact\footnote{i.e., there are no compact subsets aside from finite subsets.},  and as a result there is no hope to obtain meaningful analytic results in locally $L^0(\F)$-convex modules if one relies on classical compactness.
 A suitable substitute of compactness in modules is the so-called \emph{stable compactness}, which allows for different module analogues of classical compactness results such as a Banach-Alaouglu type theorem; see~\cite{jamzap}. 
  Shortly speaking, a subset $S\subset L^{\phi^\ast}_\F(\EE)$ is stably compact is every filter base which is closed under mixtures along $\F$-measurable countable partitions has a cluster point $y\in S$. 
 We establish a simple criterion for stable compactness in $L^{\phi^\ast}_{\F}(\EE)$.  Specifically, we prove that $S\subset L^{\phi^\ast}_{\F}(\EE)$ is stably compact if $\esssup_{y\in S}\mathbb{E}[|x_ny||\mathcal{F}] \downarrow 0$ a.s. whenever $x_n\downarrow 0$ a.s. in $L^{\phi}_{\F}(\EE)$.

 Finally, we apply the latter stable compactness criterion to obtain the robust representation of a conditional risk measure $\rho$ on the Orlicz space $L^\phi(\EE)$ with values in $L^1(\F)$. 
  We remark that Orlicz spaces provide a general framework of Banach lattices for applications in mathematical finance, described for instance in~\cite{biagini,biagini2,cheriditio,gao,orihuela}.  
 We show that every conditional risk measure $\rho\colon L^\phi(\EE)\to L^1(\F)$ uniquely extends to an $L^0(\F)$-convex function $\bar{\rho}\colon L^\phi_{\F}(\EE)\to L^0(\F)$. 
 Moreover, we extend the so-called Jouini-Schachermayer-Touzi theorem \cite{schachermayer} by showing that the dual  representation 
\begin{equation}
\rho(x)=\underset{y\in L^{\phi^\ast}(\EE)}\esssup\left\{ \mathbb{E}[x y|\F]-\rho^\#(y)\right\}
\end{equation} 
 of a conditional risk measure $\rho\colon L^\phi(\EE)\to L^1(\F)$ is attained for every $x\in L^\phi(\EE)$ if and only if $\rho$ has the Lebesgue property, and if and only if the sublever set $\{\bar{\rho}^\#\le\eta\}$ is stably compact for all $\eta\in L^0(\F)$, where $\bar{\rho}^\#$ is the   conjugate of $\bar{\rho}$ with respect to the dual system of modules $\langle L^{\phi^\ast}_\F(\EE),L^{\phi}_\F(\EE) \rangle$. 
 The proof of this result is obtained by combining the stable compactness criterion described above with scalarization techniques similar to those in~\cite{DetlefsenScandolo}.

The paper is structured as follows: Section \ref{sec0} presents preliminary concepts and notations. Section \ref{sec1} delves into the relationship between classical Orlicz spaces and Orlicz-type modules. In Section \ref{sec2}, we introduce the K\"{o}the dual of a Orlicz-type module and show that it can be canonically identified with the set of all  $\sigma$-order continuous module  homomorphisms. Section \ref{sec3} investigates sufficient conditions for stable weak compactness. Finally in Section \ref{sec4}, we apply the previous theory to obtain a robust representation result for conditional risk measures on Orlicz spaces.

\section{Preliminaries}\label{sec0}
Throughout, let $(\Omega,\mathcal{E},\mathbb{P})$ be a probability space.   
Random variables and sets which coincide $\mathbb{P}$-almost surely (a.s.) are identified. 
All inequalities between $\R\cup\{\pm\infty\}$-valued random variables are understood in the a.s.~sense. 
In the sequel, we fix a sub-$\sigma$-algebra $\mathcal{F}$ of $\mathcal{E}$, and denote by $L^0(\mathcal{F})$, $L^0(\F;\N)$, $L^0(\F;\Q)$, $L^0_+(\mathcal{F})$,  $L^0_{++}(\mathcal{F})$ and $\bar{L}^0(\mathcal{F})$ the spaces of $\P$-almost surely equal $\mathcal{F}$-measurable random variables with values in $\mathbb{R}$, $\N$, $\Q$, $[0,+\infty)$, $(0,+\infty)$, and $\R\cup\{\pm\infty\}$, respectively. 
Recall that every non-empty subset $S$ of $\bar{L}^0(\mathcal{F})$ has a least upper bound  $\esssup S$ and a greatest lower bound $\essinf S$ in $\bar{L}^0(\mathcal{F})$  with respect to the a.s. order. 
In particular, $L^0(\F)$ equipped with this order is a Dedekind complete vector lattice.  
The classical conditional expectation $\mathbb{E}[\cdot|\mathcal{F}]$, which is  defined for random variables with finite expectation, is extended to the \emph{generalized conditional expectation} $\mathbb{E}[\cdot|\mathcal{F}]\colon L^0(\mathcal{E})\to \bar{L}^0(\mathcal{F})$ by
\[
\mathbb{E}[\eta|\mathcal{F}]:=\lim_{n\to+\infty} \mathbb{E}[\eta^+ \wedge n|\mathcal{F}] - \lim_{n\to+\infty} \mathbb{E}[\eta^- \wedge n|\mathcal{F}]  
\]
with the convention $+\infty-\infty=+\infty$. 

We say that an $L^0(\F)$-module $E$ is $\F$-stable if for every sequence $(x_k)_{k\in\N}\subset E$ and every countable partition $(A_k)_{k\in\N}\subset\F$ of $\Omega$ there exists a unique $x\in E$, denoted by $\sum_{k\in\N}1_{A_k}x_k$, such that $1_{A_k}x=1_{A_k}x_k$ for all $k\in\N$. 
We call $\sum_{k\in\N}1_{A_k}x_k$ the \emph{mixture} of the sequence $(x_k)_{k\in\N}$ along the partition $(A_k)_{k\in\N}$.  
Given a non-empty subset $S$ of a $\F$-stable $L^0(\F)$-module $E$, we define the $\mathcal{F}$-\emph{stable hull} of $S$ as
\[
s_{\mathcal{F}}(S):=\left\{ \sum_{k\in\mathbb{N}} 1_{A_k}x_k\colon (x_k)_{k\in\mathbb{N}}\subset S,\:(A_k)_{k\in\mathbb{N}}\subset\F\mbox{ is a countable partition of }\Omega\right\}.
\]
The subset $S$ is said to be $\F$-\emph{stable}, if $S=s_{\mathcal{F}}(S)$. 

We have that $L^0(\EE)$ is an $\F$-stable $L^0(\F)$-module. 
In fact, given a sequence $(x_k)_{k\in\N}\subset L^0(\EE)$ and a countable partition $(A_k)_{k\in\N}\subset \F$ of $\Omega$,  the corresponding mixture is given by $\sum_{k\in\N} 1_{A_k}x_k=\lim_{n\to+\infty}\sum_{k=1}^n 1_{A_k}x_k$~a.s..

\section{Orlicz spaces and Orlicz-type modules}\label{sec1}
Let $\phi\colon[0,+\infty)\to[0,+\infty]$ be a Young function (i.e. an increasing left-continuous convex function finite on a neighborhood of $0$ with $\phi(0) = 0$ and $\lim_{x\to+\infty}\phi(x)=+\infty$), which will remain fixed throughout this note.   
Recall that the Orlicz space of $\mathcal{E}$-measurable random variables  associated with  $\phi$ is defined as
$$
L^{\phi}(\mathcal{E}):=\left\{x\in L^0(\mathcal{E})\colon \exists r\in(0,+\infty),\:\mathbb{E}[\phi(r|x|)]<+\infty \right\}.
$$
Then, $L^{\phi}(\mathcal{E})$ is a \emph{solid} vector subspace (order ideal) of $L^0(\mathcal{E})$, that is, a vector subspace such that $x\in L^{\phi}(\mathcal{E})$ whenever $|x|\le|y|$ and $y\in L^{\phi}(\mathcal{E})$. 
In addition, it is satisfied that $L^\infty(\mathcal{E})\subset L^{\phi}(\mathcal{E})\subset L^1(\mathcal{E})$. 
Both $L^p(\mathcal{E})$ ($1\le p<+\infty$) and $L^\infty(\mathcal{E})$ are particular instances of Orlicz spaces. 
Given the sub-$\sigma$-algebra $\mathcal{F}$ of $\mathcal{E}$, the \emph{Orlicz-type} $L^0(\F)$-module  associated with  $\phi$  is defined as
\begin{equation}
\label{eq:OrliczModule}
L^{\phi}_\F(\mathcal{E}):=\left\{x\in L^0(\mathcal{E})\colon \exists \eta\in L^0_{++}(\F),\:\mathbb{E}[\phi(\eta|x|)|\F]\in L^0(\F) \right\}.
\end{equation}
Then, $L^{\phi}_\F(\mathcal{E})$ is both an $L^0(\F)$-submodule of $L^0(\mathcal{E})$ and an order ideal of $L^0(\mathcal{E})$. 
In addition, Vogelpoth~\cite{vogelpoth2009phd} proved the following  product structure that relates $L^\phi(\mathcal{E})$ and $L^\phi_{\F}(\mathcal{E})$ given by
\begin{equation}
\label{eq:productOrlicz}
L^\phi_\F(\mathcal{E})=L^0(\F)L^\phi(\mathcal{E}).
\end{equation}

The following result gives a further relation between $L^\phi(\mathcal{E})$ and $L^\phi_{\F}(\mathcal{E})$. 

\begin{proposition}\label{prop: ccHull}
 $L^\phi_\mathcal{F}(\mathcal{E})$ is $\F$-stable. 
Moreover,
 $s_\F\left(L^\phi(\mathcal{E})\right)=L^\phi_\mathcal{F}(\mathcal{E}).$
\end{proposition}

\begin{proof} 
Fix a sequence $(x_k)_{k\in\N}\subset L^\phi_{\F}(\EE)$ and a countable partition $(A_k)_{k\in\N}\subset\F$ of $\Omega$. 
Let us prove that $x:=\sum_{k\in\N} 1_{A_k}x_k\in L^\phi_{\F}(\EE)$. 
It follows from \eqref{eq:OrliczModule} that for every $k\in\N$ there exists $\eta_k\in L^0_{++}(\F)$ such that $\mathbb{E}[\phi(\eta_k|x_k|)|\F]\in L^0(\F)$. 
Now, let  $\eta:=\sum_{k\in\N}1_{A_k}\eta_k\in L^0_{++}(\F)$.  We have
\begin{align*}
\mathbb{E}[\phi(\eta|x|)|\F]&=\sum_{k\in\N} 1_{A_k}\mathbb{E}[\phi(\eta|x|)|\F]\\
&=\sum_{k\in\N} 1_{A_k}\mathbb{E}[1_{A_k}\phi(\eta|x|)|\F]\\
&=\sum_{k\in\N} 1_{A_k}\mathbb{E}[1_{A_k}\phi(\eta_k|x_k|)|\F]\in L^0(\F).
\end{align*}
This implies that $x\in L^\phi_{\F}(\EE)$  in view of \eqref{eq:OrliczModule}.  
We conclude that  
$L^\phi_\mathcal{F}(\mathcal{E})$ is $\F$-stable. 
Since  $L^\phi(\EE)\subset L^\phi_\F(\EE)$, one has that   
$s_\F\left(L^\phi(\EE)\right)\subset s_\F(L^\phi_\F(\EE))=L^\phi_\F(\EE)$. 
As for the other inclusion, suppose that $x\in L^\phi_\F(\EE)$. 
Due to \eqref{eq:productOrlicz}, we have that $x=\eta x_0$ for certain $\eta\in L^0(\F)$ and $x_0\in L^\phi(\EE)$. 
We have that $x=\sum_{k\in\N} 1_{A_k}\eta x_0$, where $A_k:=\{k-1\leq|\eta|<k\}$ for each $k\in\N$. 
 Since $|1_{A_k}\eta x_0|\leq k |x_0|$ and $L^\phi(\EE)$ is an order ideal of $L^0(\mathcal{E})$, it follows that $x_k:=1_{A_k}\eta x_0\in L^\phi(\EE)$. 
 We conclude that $x=\sum_{k\in\N}1_{A_k}x_k\in s_\F\left(L^\phi(\EE)\right)$. 
\end{proof}

\section{Duality}\label{sec2}
The conjugate Young function of $\phi$ is defined as 
$$\phi^\ast\colon [0,+\infty)\to [0,+\infty],\quad\phi^\ast(s):=\sup_{r\ge 0}\{r s - \phi(r)\}.$$  
It is known that $L^{\phi^\ast}(\EE)$ agrees with the  \emph{K\"{o}the dual} of $L^{\phi}(\EE)$; see, e.g.~\cite{nowak,owari2014lebesgue,zaanenII}. 
Namely, the following holds true
\begin{equation}
\label{eq:KotheDual}
L^{\phi^\ast}(\EE)=\left\{ y\in L^0(\mathcal{E})  \colon \E[|x y|]<+\infty\text{ for all }x\in  L^{\phi}(\EE) \right\}.
\end{equation}

In the following, we extend \eqref{eq:KotheDual} to Orlicz-type modules. 
We define the dual K\"{o}the module of $L^{\phi}_{\F}(\EE)$  as 
\[
\big(L^{\phi}_{\F}(\EE) \big)^\#:=\left\{ y\in L^0(\mathcal{E})  \colon \E[|x y||\F]<+\infty\mbox{ for all }x\in  L^{\phi}_{\F}(\EE) \right\}.
\] 
Then we have the following.
\begin{proposition}\label{prop:dualKothe}
$
\big(L^{\phi}_{\F}(\EE) \big)^\#=L^{\phi^\ast}_{\F}(\EE).
$
\end{proposition}

To prove Proposition~\ref{prop:dualKothe}, we need some preliminary results. 
Suppose that $E$ is an $L^0(\F)$-module.   
A function $\Vert\cdot\Vert\colon E\to L^0_+(\F)$ is called an $L^0(\F)$-\emph{norm} if:
\begin{itemize}
\item[(i)] $\Vert  x+y\Vert\le \Vert x\Vert +\Vert y\Vert$,
\item[(ii)] $\Vert \eta x\Vert=|\eta|\Vert x\Vert$ for all~$\eta\in L^0(\F)$,
\item[(iii)] $\Vert x\Vert=0$ implies $x=0$.
\end{itemize}
The collection of all balls
 $B_\varepsilon(x):=\left\{y\in E \colon \Vert y-x\Vert\le \varepsilon \right\}$ 
such that $\varepsilon\in L^0_{++}(\F)$ and $x\in E$ is a base for a topology on $E$, which we call the topology induced by $\Vert\cdot\Vert$; see~\cite{kupper03}.  
The function $\eta\mapsto|\eta|$ is an $L^0(\F)$-norm on $L^0(\F)$.  
Hereafter, we assume that $L^0(\F)$ is endowed with the topology induced by $|\cdot|$.  

Given an $L^0(\F)$-normed module $(E,\Vert\cdot\Vert)$, we define its \emph{topological dual}, denoted by $E^\ast$, as the collection of all  module homomorphisms $\mu\colon E\to L^0(\F)$ which are $L^0(\F)$-norm continuous.  
The space $E^\ast$ is an $L^0(\F)$-normed module. 
More specifically, for $\mu\in E^\ast$ we set
\begin{equation}\label{eq:dualNorm}
\Vert\mu\Vert:=\esssup\left\{|\mu(x)| \colon \Vert x\Vert\le 1,\: x\in E\right\},
\end{equation}
which defines an $L^0(\F)$-norm on $E$. 
The $L^0(\F)$-module $L^{\phi}_{\F}(\EE)$ can be endowed with the following $L^0(\F)$-norm
\[
\Vert x|\F\Vert_\phi:=\essinf\{\lambda\in L^0_{++}(\F)\colon \mathbb{E}[\phi(\tfrac{|x|}{\lambda})|\F]\le 1\};
\]
 see~\cite{vogelpoth2009phd}, which generalizes the classical  Luxemburg norm.  
 In particular, we consider the topological dual $\big(L^{\phi}_{\F}(\EE)\big)^\ast$, which is endowed with an  $L^0(\F)$-norm defined as in~\eqref{eq:dualNorm}. 
 
\begin{definition}
Let  $f\colon L^{\phi}_{\F}(\EE)\to \bar{L}^0(\F)$ be a function.
\begin{enumerate}
\item $f$ has the \emph{local property} if $1_A f(x)=1_A f(1_A x)$ for all $x\in L^{\phi}_{\F}(\EE)$ and $A\in\F$.
\item $f$ is increasing if $f(x)\le f(y)$ whenever $x\le y$. 
\item $f$ is $L^0(\F)$-convex if $f(\lambda x + (1-\lambda)y)\le \lambda f(x)+(1-\lambda)f(y)$ for all $x,y\in L^{\phi}_{\F}(\EE)$ and $\lambda\in L^0(\F)$ with $0\le \lambda\le 1$.
\end{enumerate}
\end{definition} 

The $L^0(\F)$-module $L^{\phi}_{\F}(\EE)$ is complete in the sense that every Cauchy net in $L^{\phi}_{\F}(\EE)$ is convergent; see~\cite[Theorem 3.3.3]{vogelpoth2009phd}.  
The following result was proven in~\cite[Theorem 4.1.3]{vogelpoth2009phd} for general complete  $L^0(\F)$-normed module lattices.  
Next, we state the result for the complete module $L^{\phi}_{\F}(\EE)$.
 \begin{lemma}\label{lem:autCont}
 Let $f\colon L^{\phi}_{\F}(\EE)\to L^0(\F)$ be a function with the local property. 
 If $f$ is increasing and convex, then $f$ is $L^0(\F)$-norm continuous. 
 \end{lemma}

\begin{lemma}\label{lem:autCont2}
Let $f\colon L^{\phi}_{\F}(\EE)\to L^0(\F)$ be a function which has the local property and is increasing.  
 If $f$ is convex, then $f$ is $L^0(\F)$-convex. 
 If $f$ is linear, then $f$ is a module homomorphism. 
\end{lemma}
\begin{proof}
We proceed in two steps. 

\emph{Step 1}: Fix $\lambda\in L^0(\F)$ and  any open neighborhood $U$ of $\lambda$ in $(L^0(\F),|\cdot|)$. 
In this step we prove that there exists  $\lambda_U\in L^0(\F;\Q)$ such that $\lambda_U\in U$.  
Moreover, if $0\le\lambda\le 1$ we can  choose $\lambda_U$ so that $0\le \lambda_U\le 1$. 
 To prove that take $\varepsilon\in L^0_{++}(\F)$ such that $B_{\varepsilon}(\lambda)\subset U$. 
  Given a enumeration $\mathbb{Q}:=\{q_1,q_2,\cdots\}$  of the set of rational numbers, we define $A_1:=\{\lambda\le q_1<\lambda + \varepsilon\}$ and  $A_k:=\{\lambda\le q_k<\lambda +\varepsilon\}\setminus\bigcup_{i=1}^{k-1}A_i$ for $k=2,3,\cdots$.  
  Due to the density of $\Q$, we have that $(A_k)_{k\in\N}$ is a partition of $\Omega$. 
  We define $\lambda_U:=\sum_{k\in\N} 1_{A_k} q_k$. 
  By construction, we have that $\lambda_U\in B_{\varepsilon}(\lambda)\subset U$ as desired. 
  Finally, note that if $0\le \lambda\le 1$, then $\lambda_U\wedge 1$ meets the required conditions. 
 
\emph{Step 2}: Let $\lambda\in L^0(\F;\Q)$ and suppose that $f$ has the local property. 
\begin{enumerate}
\item If $f$ is convex and $0\le \lambda\le 1$, then $f(\lambda x + (1-\lambda)y)\le \lambda f(x)+(1-\lambda)f(y)$ for all $x,y\in L^\phi_\F(\EE)$. 
In fact, fixed $q\in\mathbb{Q}$, due to the local property
\begin{align*}
1_{\{\lambda=q\}}f(\lambda x + (1-\lambda)y)&=1_{\{\lambda=q\}}f(1_{\{\lambda=q\}}(\lambda x + (1-\lambda)y))\\
&=1_{\{\lambda=q\}}f(1_{\{\lambda=q\}}(q x + (1-q)y))\\
&=1_{\{\lambda=q\}}f(q x + (1-q)y)\\
&\le 1_{\{\lambda=q\}}\Big(q f(x) + (1-q)f(y)\Big)\\
&= 1_{\{\lambda=q\}}\Big(\lambda f(x) + (1-\lambda)f(y)\Big).
\end{align*}
Summing up over all $q\in \mathbb{Q}$ we obtain the desired conclusion.	
\item If $f$ is linear, then $f(\lambda x)=\lambda f(x)$. 
This is proven similarly as in (i). 
\end{enumerate}
Finally, since $f$ is increasing, we have that $f$ is $L^0(\F)$-norm continuous by Lemma~\ref{lem:autCont}. 
 We conclude that (i) and (ii) above are valid for every $\lambda\in L^0(\F)$ due to the step 1 and the continuity of $f$.  
 The proof is complete. 
\end{proof}

 For each $y\in \big(L^{\phi}_{\F}(\EE) \big)^\#$, we define the function
 \[
 \mu_y\colon L^{\phi}_{\F}(\EE)\to L^0(\F),\quad \mu_y(x):=\mathbb{E}[x y|\F].
 \]
\begin{lemma}\label{lem:ExpnormCont}
If $y\in \big(L^{\phi}_{\F}(\EE) \big)^\#$, then $\mu_y\in \big(L^{\phi}_{\F}(\EE) \big)^\ast$. 
In particular, 
\[
|\mathbb{E}[x y|\F]|\le \Vert\mu_y\Vert \Vert y|\F\Vert_\phi\quad\mbox{ for all }x\in L^{\phi}_{\F}(\EE).
\]
\end{lemma}
\begin{proof}
Suppose first that $y\ge 0$. 
Then, $\mu_y$ is increasing and convex. 
Due to Lemma~\ref{lem:autCont}, $\mu_y$ is an $L^0(\F)$-norm continuous module homomorphism. 
For arbitrary $y$, we put $y=y^+ - y^-$ and apply the previous case to $y^+$ and $y^-$. 
The proof is complete. 
\end{proof}

\begin{lemma}\label{lem:finExp}
If $x\in L^\phi(\EE)$, then $\mathbb{E}[ \Vert x|\F\Vert_\phi ]<+\infty.$
\end{lemma}
\begin{proof}
Let $r\in(0,+\infty)$ satisfy that $\mathbb{E}[\phi(|x|r^{-1})]<+\infty$. 
By convexity of $\phi$, 
\[
\mathbb{E}\left[\phi\left(\frac{|x|r^{-1}}{1+\mathbb{E}[\phi(|x|r^{-1})|\F]}\right)\Big| \F \right]
\le \frac{\mathbb{E}\left[\phi\left(|x|r^{-1}\right)|\F\right]}{1+\mathbb{E}[\phi(|x|r^{-1})|\F]}\le 1.
\]
 It follows that
 \[
 \frac{1}{r}\Vert x|\F\Vert_\phi\le 1+\mathbb{E}[\phi(|x|r^{-1})|\F].
 \]
Then, taking expectations
 \begin{align*}
 \frac{1}{r}\mathbb{E}\Big[\Vert x|\F\Vert_\phi\Big]&\le 1+\mathbb{E}\Big[\mathbb{E}[\phi(|x|r^{-1})|\F]\Big]\\
 &=1+\mathbb{E}\Big[\phi(|x|r^{-1})\Big]<+\infty.
\end{align*}
The proof is complete.
\end{proof}

Now, we turn to the proof of Proposition~\ref{prop:dualKothe}.

\begin{proof}
Due to Proposition~\ref{prop: ccHull}, we have that 
$s_{\F}(L^{\phi^\ast}(\EE))=L^{\phi^\ast}_{\F}(\EE)$. 
Then, it suffices to prove that
\begin{equation}\label{eq:stableDual}
s_{\F}(L^{\phi^\ast}(\EE))=\big(L^{\phi}_{\F}(\EE) \big)^\#.
\end{equation}
It follows by inspection that $\big(L^{\phi}_{\F}(\EE) \big)^\#$ is $\F$-stable. 
Hence, to prove `$\subset$' in \eqref{eq:stableDual} it is enough to show that  $L^{\phi^\ast}(\EE)\subset \big(L^{\phi}_{\F}(\EE) \big)^\#$. 
Fix $y\in L^{\phi^\ast}(\EE)$. 
If $x\in L^{\phi}_{\F}(\EE)$, by Proposition~\ref{prop: ccHull}, we have that $x=\sum_{k\in\N} 1_{A_k} x_k$ for a sequence $(x_k)_{k\in\N}\subset L^{\phi}(\EE)$ and a countable partition $(A_k)_{k\in\N}\subset \F$ of $\Omega$. 
Therefore, $xy=\sum_{k\in\N} 1_{A_k} x_k y$, and
\[
\mathbb{E}[|x y||\F]=\sum_{k\in\N} 1_{A_k} \mathbb{E}[x_k y|\F]<+\infty,
\] 
where the right hand side is finite by~\eqref{eq:KotheDual}. 
Since $x$ was arbitrary, we conclude that $y\in \big(L^{\phi}_{\F}(\EE) \big)^\#$.  

As for the other inclusion, suppose that $y\in \big(L^{\phi}_{\F}(\EE) \big)^\#$. 
Due to Lemma~\ref{lem:ExpnormCont}, we know that $\Vert\mu_y\Vert<+\infty$. 
We can find a sequence $(r_k)_{k\in\N}$ of positive real numbers and a countable partition $(A_k)_{k\in\N}\subset\F$ of $\Omega$ such that $\Vert\mu_y\Vert<\sum_{k\in\N} 1_{A_k}r_k$. 
For each $k\in\N$, define $y_k:=1_{A_k}y$. 
We prove that $y_k\in L^{\phi^\ast}(\EE)$ for all $k$. 
In view of \eqref{eq:KotheDual}, it suffices to show that $\mathbb{E}[|x y_k|]<+\infty$ for all $x\in L^\phi(\EE)$. 
 Indeed, if $x\in L^\phi(\EE)$, we have
 \begin{align*}
 \mathbb{E}[|x y_k||\F]&=1_{A_k}\mathbb{E}[1_{A_k} |x|{\rm sgn}(y)y|\F]\\
 &=1_{A_k}\mu_y(|x|1_{A_k}{\rm sgn}(y))\\
 &\le 1_{A_k}\Vert\mu_y\Vert \Vert 1_{A_k}|x||\F\Vert_\phi\\
  &\le 1_{A_k} r_k  \Vert x|\F\Vert_\phi\\
  &\le r_k  \Vert x|\F\Vert_\phi.  
 \end{align*}
Taking expectations, we get
\[
\mathbb{E}[|x y_k|]\le r_k \mathbb{E}\left[\Vert x|\F\Vert_\phi\right]<+\infty,
\]
where  the right hand side is finite due to Lemma~\ref{lem:finExp}. The proof is complete.
\end{proof}

\begin{definition}
Let $\mathcal{X}\subset L^0(\EE)$ be a vector sublattice and $\mu\colon \mathcal{X}\to L^0(\F)$ a linear function.   
We say that $\mu$ is \emph{order bounded} if it maps order bounded sets onto order bounded sets.    
We say that $\mu$  is \emph{$\sigma$-order continuous} if $\essinf|\mu(x_n)|=0$ whenever $x_n\downarrow 0$ a.s. in $\mathcal{X}$. 
\end{definition}

\begin{proposition}\label{prop:OrdContNormCont}
Let $\mu\colon L^\phi_{\F}(\mathcal{E})\to L^0(\F)$ be a module homomorphism. 
Then, $\mu$ is order bounded if and only if $\mu \in \big(L^{\phi}_\F(\mathcal{E})\big)^\ast$. 
\end{proposition}
\begin{proof}
If $\mu\in \big(L^{\phi}_\F(\mathcal{E})\big)^\ast$, we have that $|\mu(x)|\le \Vert\mu\Vert\Vert x|\F\Vert_\phi$ for all $x\in L^\phi_\F(\EE)$, which implies that $\mu$ is order bounded. 
Conversely, suppose that $\mu$ is order bounded. 
In that case, since $L^0(\F)$ is Dedekind complete, then $\mu=\mu^+ - \mu^-$ where $\mu^+,\mu^- \colon L^\phi_{\F}(\mathcal{E})\to L^0(\F)$ are positive linear functions; see \cite[Theorem 20.2]{zaanen}. 
More precisely, $\mu^+(x)=\esssup\{\mu(y) \colon 0\le y\le x\}$ for $x\ge 0$, and for arbitrary $x\in L^\phi_{\F}(\mathcal{E})$  we have $\mu^+(x)=\mu^+(x^+)-\mu^+(x^-)$ and $\mu^-(x)=(-\mu)^+(x)$; see \cite[Theorem 20.4]{zaanen}. 
Both $\mu^+$ and $\mu^-$ have the local property, and are linear and increasing. 
We conclude that $\mu^+,\mu^-\in \big(L^{\phi}_\F(\mathcal{E})\big)^\ast$ due to Lemma~\ref{lem:autCont} and  Lemma~\ref{lem:autCont2}.    
Consequently, $\mu\in \big(L^{\phi}_\F(\mathcal{E})\big)^\ast$. 
\end{proof}

The following lemma is well-known in the theory of Orlicz spaces; see~e.g.~\cite[Section 3]{nowak} and~\cite{owari2014lebesgue,zaanenII}.
\begin{lemma}\label{lem:repLin}
A linear function $\mu\colon L^\phi(\mathcal{E})\to\R$ is $\sigma$-order continuous  if and only if there exists (a unique) $y\in L^{\phi^\ast}(\EE)$ such that $\mu(x)=\mathbb{E}[x y]$ for all $x\in L^\phi(\mathcal{E})$.
\end{lemma}

The next result extends Lemma~\ref{lem:repLin} to the present setting.
\begin{theorem}\label{thm:ordCont}
Let $\mu\colon L^\phi_{\F}(\mathcal{E})\to L^0(\F)$ be a module homomorphism. 
Then, $\mu$ is $\sigma$-order continuous and order bounded if  and only if there exists (a unique) $y\in L^{\phi^\ast}_\F(\mathcal{E})$ such that $\mu=\mu_y$.
\end{theorem}
\begin{proof}
Suppose that there exists $y\in L^{\phi^\ast}_\F(\mathcal{E})$ such that $\mu=\mu_y$. 
Then, $\mu$ is order bounded and, moreover, it follows by monotone convergence that $\mu$ is $\sigma$-order continuous. 
Conversely, suppose that $\mu$ is $\sigma$-order continuous and order bounded.  
By Proposition \ref{prop:OrdContNormCont}, we have that $\Vert \mu\Vert<+\infty$. 
For each $n\in\N$, define $A_n:=\{n-1 \le \Vert\mu\Vert<n\}$.  
For every $n\in\N$ and $x\in L^{\phi}(\mathcal{E})$, we have that 
\[
\E[1_{A_n}\mu(x)]\leq n\E[\Vert x|\F\Vert_\phi]<+\infty,
\]
where the right hand side is finite due to Lemma~\ref{lem:finExp}. 
Therefore, the function
\[
\mu_n\colon L^\phi(\mathcal{E})\rightarrow \R,\quad \mu_n(x):=\E[1_{A_n}\mu(x)],
\]
is well-defined for all $n\in\N$. 
Moreover, $\mu$ is linear and order bounded, and, by monotone convergence, the $\sigma$-order continuity of $\mu$ implies that $\mu_n$ is $\sigma$-order continuous. 
Then, by Lemma \ref{lem:repLin}, for each $n\in\N$ there exists $y_n\in L^{\phi^\ast}(\mathcal{E})$ such that $\mu_n(x)=\E[x y_n]$ for all $x\in L^\phi(\mathcal{E})$. 
In addition, for each $n\in\N$ we have that
\begin{equation}\label{eq:homRep}
1_{A_n}\mu(x)=\mathbb{E}[1_{A_n}  x y_n|\F]\quad\mbox{ for all }x\in L^{\phi}(\EE).
\end{equation}
In fact, for each $A\in\F$ one has 
$$\mathbb{E}[1_A 1_{A_n}\mu(x)]=\mathbb{E}[ 1_{A_n}\mu(1_A 1_{A_n}x)]=\mu_n(1_A 1_{A_n} x)=\mathbb{E}[1_A 1_{A_n} x y_n],$$
 then \eqref{eq:homRep} follows by definition of conditional expectation.  
If we set $y=\sum_{k\in\N} 1_{A_k} y_k$, then, due to \eqref{eq:homRep}, for every $x\in L^\phi(\mathcal{E})$ we have
\begin{align*}
\mu(x)&=\sum_{k\in\N} 1_{A_k}\mu(x)\\
&=\sum_{k\in\N} 1_{A_k}\mathbb{E}[1_{A_{k}} x y|\F]\\
&=\sum_{k\in\N} 1_{A_k}\mathbb{E}[x y|\F]\\
&=\mathbb{E}[x y|\F].
\end{align*}
 Proposition~\ref{prop: ccHull} allows to extend the latter to any $x\in L^\phi_{\F}(\mathcal{E})$. 
Consequently, $\mu=\mu_y$. 
Finally, if $\mu=\mu_y=\mu_z$, then 
\[
\mathbb{E}[|y-z||\F]\le \left|\mathbb{E}[y-z|\F]\right|=\left|\mathbb{E}[y|\F] -\mathbb{E}[z|\F]\right|=|\mu_y(1)-\mu_{z}(1)|=0. 
\] 
Then, $\mathbb{E}[|y-z||\F]=0$ and, in turn, $y=z$. 
\end{proof}

\section{A criterion for stable compactness}\label{sec3}
In the following, we recall the concepts of stable weak topology and  stable weak compactness; see \cite{jamzap}. 
Suppose that we have a bi-$L^0(\F)$-linear form $\langle\cdot,\cdot\rangle\colon E\times F\to L^0(\F)$ on two $\F$-stable $L^0(\F)$-modules $E$ and $F$.\footnote{i.e., the mappings $x\mapsto \langle x,y_0\rangle$ and $y\mapsto \langle x_0,y\rangle$ are module homomorphisms for all $x_0\in E$ and $y_0\in F$. Recall that an $L^0(\F)$-module $E$ is stable if for every sequence $(x_k)\subset E$ and every countable partition $(A_k)\subset \F$ of $\Omega$, there exists a unique $x\in E$, denoted by $\sum_{k\in\N}1_{A_k}x_k$, such that $1_{A_k}x=1_{A_k}x_k$ for all $k\in\N$.}  
For every sequence $(F_k)_{k\in\N}$ of non-empty finite subsets of $F$, every countable partition $(A_k)_{k\in\N}\subset\F$ of $\Omega$, and every $\varepsilon\in L^0_{++}(\F)$ define
\begin{equation}\label{eq:basicNeigh}
U_{(F_k),(A_k),\varepsilon}:=\left\{ x\in E \colon \sum_{k\in\N} 1_{A_k}\underset{y\in F_k}\esssup|\langle x,y\rangle|<\varepsilon \right\}.
\end{equation}
The collection of sets 
$$\left\{ x+U_{(F_k),(A_k),\varepsilon}\colon
\begin{array}{c}
x\in E,\: \emptyset\neq F_k\subset F\mbox{ is finite, }\varepsilon\in L^0_{++}(\F)\\
(A_k)\subset\F\mbox{ is a countable partition of }\Omega\\
\end{array}
\right\}$$
is a base for a topology on $E$ which is called the stable weak topology on $E$ induced by $\langle\cdot,\cdot\rangle$, and is denoted by $\sigma_s(E,F)$.

\begin{remark}\label{rem:stTop}
The stable weak topologies were studied in \cite{orihuela3,jamzap}. 
In particular, they are stable locally $L^0(\F)$-convex topologies as defined in \cite{orihuela3,jamzap}.  
A fundamental property of stable locally $L^0(\F)$-convex topologies is that, for every countable partition $(A_k)_{k\in\N}\subset\F$ of $\Omega$ and every sequence $(O_k)_{k\in\N}$ of non-empty open sets,  the mixture $\sum_{k\in\N}1_{A_k}O_k$ is again an open set; see~\cite[Theorem 1.1]{orihuela3}.
\end{remark}
 
Suppose that $E$ is endowed with the stable weak topology $\sigma_s(E,F)$. We introduce the following terminology: 

\begin{itemize}
\item An $\F$-\emph{stable filter base} on $E$ is a non-empty collection $\mathcal{B}$ of non-empty subsets of $E$ such that: 
\begin{itemize}
\item[(i)] $\mathcal{B}$ is a filter base, i.e.~for any $U,V\in\mathcal{B}$, there exists $W\in\mathcal{B}$ such that $W\subset U\cap V$;
\item[(ii)] each $U\in\mathcal{B}$ is $\F$-stable;
\item[(iii)] if $(U_k)_{k\in\N}\subset \mathcal{B}$, and $(A_k)_{k\in\N}\subset\F$ is a countable partition of $\Omega$, then $\sum_{k\in\N} 1_{A_k}U_k$ 
is an element of $\mathcal{B}$.\footnote{We define $\sum_{k\in\N} 1_{A_k}U_k:=\left\{\sum_{k\in\N} 1_{A_k} y_k\colon y_k\in U_k\mbox{ for all }k\in\N\right\}$.}
\end{itemize}
\item We say that $x\in E$ is a cluster point of the $\F$-stable filter base $\mathcal{B}$ on $E$ if for every $U\in \mathcal{B}$ it holds $x\in {\rm cl}(U)$.\footnote{${\rm cl}(U)$ denotes the topological closure of $U$.}
\item
An $\F$-stable subset $S$ of $E$ is said to be \emph{stably compact} if any $\F$-stable filter base $\mathcal{B}$ consisting of subsets of $S$ has a cluster point $x\in S$. 
\item An $\F$-stable subset $S$ of $E$ is said to be \emph{relatively stably compact} if ${\rm cl}(S)$ is stably compact. 
\end{itemize}

In the sequel, we focus on three particular instances of stable weak topologies:
\begin{itemize}
\item The $w_s$-topology on $L^\phi_\F(\EE)$, which is induced by the bi-$L^0(\F)$-linear form  $$(x,y)\mapsto \mathbb{E}[xy|\F]\colon L^\phi_\F(\EE)\times L^{\phi^\ast}_\F(\EE)\to L^0(\F).$$
\item  The $w_s^\#$-topology on $L^{\phi^\ast}_\F(\EE)$, which is induced by the bi-$L^0(\F)$-linear form  $$(y,x)\mapsto \mathbb{E}[xy|\F]\colon L^{\phi^\ast}_\F(\EE)\times L^{\phi}_\F(\EE)\to L^0(\F).$$
\item The $w_s^\ast$-topology on $\big(L^\phi_\F(\EE)\big)^\ast$, which is induced by the bi-$L^0(\F)$-linear form  $$(\mu,x)\mapsto \mu(x)\colon \big(L^\phi_\F(\EE)\big)^\ast\times L^{\phi}_\F(\EE)\to L^0(\F).$$
\end{itemize}
 Notice that, under the canonical embedding $$y\mapsto \mu_y\colon L^{\phi^\ast}_\F(\EE)\to \big(L^\phi_\F(\EE)\big)^\ast,$$ the $w^\#_s$-topology is the restriction of the $w^\ast_s$-topology to $L^{\phi^\ast}_\F(\EE)$.

We next state the main result of this section.
\begin{theorem}\label{thm:criterion}
Let $S$ be an $\F$-stable subset of $L^{\phi^\ast}_\F(\mathcal{E})$. 
If 
\begin{equation}
\label{eq:criterion}
\underset{y\in S}\esssup\mathbb{E}\left[|x_n y||\mathcal{F}\right]\to 0\mbox{ a.s. whenever }x_n\downarrow 0\mbox{ a.s. in $L^{\phi}_\F(\mathcal{E})$,}	
\end{equation}
then $S$ is relatively stably $w_s^\#$-compact.
\end{theorem}

To prove Theorem \ref{thm:criterion} we require some preliminary results. 
 The following Banach-Alaoglu type theorem is a particular instance of~\cite[Corollary 6.8]{jamzap}.
\begin{lemma}\label{lem:alaoglu}
The set
\[
\left\{ \mu\in\big(L^\phi_{\F}(\mathcal{E})\big)^\ast \colon \Vert\mu\Vert\le 1\right\}
\]
is stably $w_s^\ast$-compact.
\end{lemma}

Next, we state the extended uniform boundedness principle \cite[Theorem 6.16]{jamzap} in the particular case of the complete $L^0(\F)$-normed module $L^{\phi}_\F(\mathcal{E})$.
\begin{lemma}\label{lem:UBP}
Let $S$ be an $\F$-stable subset of $\big(L^{\phi}_\F(\mathcal{E})\big)^\ast$ such that for every $x\in L^{\phi}_\F(\mathcal{E})$, there exists $\eta_x$ such that $\esssup_{\mu\in S}|\mu(x)|\le \eta_x$. 
Then, there exists $\eta\in L^0(\F)$ such that $\Vert \mu\Vert\le \eta$ for all $\mu\in S$.
\end{lemma}

Now we are ready to prove Theorem \ref{thm:criterion}.

\begin{proof} 
We prove that  
\begin{equation}\label{eq:closures}
{\rm cl}_{w_s^\ast}(S)={\rm cl}_{w_s^\#}(S),
\end{equation}
where ${\rm cl}_{w_s^\ast}(S)$ and ${\rm cl}_{w_s^\#}(S)$ denote the topological closures of $S$ in the $w_s^\ast$-topology and the $w_s^\#$-topology, respectively. 
Since the $w^\#_s$-topology on $L^{\phi^\ast}_\F(\mathcal{E})$ is exactly the topology induced on $L^{\phi^\ast}_\F(\mathcal{E})$ by  the $w^\ast_s$-topology,  
it is enough to show that ${\rm cl}_{{w}_s^\ast}(S)\subset L^{\phi^\ast}_\F(\mathcal{E})$.  
Fix $\mu\in {\rm cl}_{{w}_s^\ast}(S)$. 
Given $x\in L^{\phi}_\F(\mathcal{E})$ and  
 $n\in\N$ there exists $y_n\in S$ such that 
$|(\mu_{y_n}-\mu)(x)|\le 1/n$. 
Then,
\[
|\mu(x)|\le |(\mu-\mu_{y_n})(x)|+|\mu_{y_n}(x)|\le \frac{1}{n} + \underset{y\in S}\esssup\mathbb{E}\left[|x y||\mathcal{F}\right].
\]
 Letting $n\to+\infty$, we conclude that 
 \begin{equation}
 |\mu(x)|\le \esssup_{y\in S}\mathbb{E}\left[|x y||\mathcal{F}\right]\quad\mbox{ for all }x\in L^{\phi}_\F(\mathcal{E}). 
 \end{equation}
 Now, if $x_n\downarrow 0$ a.s. in $L^{\phi}_\F(\mathcal{E})$, then it follows from the inequality above and the assumption \eqref{eq:criterion} that $\esssup_{n\in\N}|\mu(x_n)|=0$.  
 This shows that $\mu$ is $\sigma$-order continuous. 
 Also, $\mu$ is order bounded because it is $L^0(\F)$-norm continuous.  
 Then, by Theorem~\ref{thm:ordCont}, we have that $\mu=\mu_y$ for certain $y\in L^{\phi^\ast}_\F(\mathcal{E})$, and \eqref{eq:closures} follows. 
 
 Fix $x\in L^\phi_\F(\EE)$. 
 Given $y\in S$, since $|\tfrac{1}{n}x|\downarrow 0$~a.s. by the assumption \eqref{eq:criterion} we have that $\tfrac{1}{n}\esssup_{y\in S}\mathbb{E}\left[|x y||\mathcal{F}\right]\to 0$. 
 In particular,  $\esssup_{y\in S}\mathbb{E}\left[|x y||\mathcal{F}\right]$ is finite. 
 By Lemma~\ref{lem:UBP}, there exists $\eta\in L^0_{++}$ such that 
 $S\subset \left\{\mu \in \big(L^\phi_{\F}(\EE)\big)^\ast\colon \Vert \mu\Vert\le\eta\right\}$. 
 By Lemma~\ref{lem:alaoglu}, we obtain that $S$ is relatively stably ${w}_s^\ast$-compact, and in turn  relatively stably ${w}_s^\#$-compact due to \eqref{eq:closures}.  
\end{proof}

\section{Application to conditional risk measures}\label{sec4}
The notion of conditional risk measure was independently introduced in
\cite{BionNadal} and \cite{DetlefsenScandolo}, where it is considered $L^\infty$ as a model space. 
Conditional risk measures on the larger space $L^p$ ($1\le p<+\infty$) were introduced in \cite{kupper11}.    
In the following, we study the dual representation of conditional risk measures on general Orlicz spaces.

\begin{definition}
A function $\rho\colon L^\phi(\mathcal{E})\to L^0(\F)$ is:
\begin{enumerate}
\item \emph{Monotone} if $x\le y$ implies $\rho(x)\ge \rho(y)$;
\item \emph{$L^\infty(\F)$-convex} if $\rho(\lambda x + (1-\lambda)y)\le \lambda\rho(x)+(1-\lambda)\rho(y)$ whenever $\lambda\in L^\infty(\F)$ and $0\le\lambda\le 1$;
\item \emph{$L^\infty(\F)$-cash invariant} if $\rho(x+\eta)=\rho(x)-\eta$ for all $\eta\in L^\infty(\F)$.
\end{enumerate}
We say that $\rho$ is a \emph{conditional risk measure} on $L^\phi(\mathcal{E})$ if $\rho$ is monotone, $L^\infty(\F)$-convex and $L^\infty(\F)$-cash invariant. 
The \emph{conjugate} of a conditional risk measure $\rho$ on $L^\phi(\mathcal{E})$ is defined as 
\begin{equation}
\label{eq:FenchelConjugate}
\rho^\#(y):=\esssup\left\{\E[x y|\F]-\rho(x) \colon x\in L^\phi(\mathcal{E})\right\}\quad\mbox{ for }y\in L^{\phi^*}(\mathcal{E}). 
\end{equation}
\end{definition}

The previous definition involves functions on vector spaces. 
 Filipovi\'{c} et al.~\cite{kupper11} also introduced a module-based approach to conditional risk measures. 
 In line with \cite{kupper11} we introduce the following terminology.

\begin{definition}
A function $\rho\colon L^\phi_{\F}(\mathcal{E})\to L^0(\F)$ is:
\begin{enumerate}
\item \emph{Monotone} if $x\le y$ implies $\rho(x)\ge \rho(y)$;
\item \emph{$L^0(\F)$-convex} if $\rho(\lambda x + (1-\lambda)y)\le \lambda\rho(x)+(1-\lambda)\rho(y)$ whenever $\lambda\in L^0(\F)$ and $0\le\lambda\le 1$;
\item \emph{$L^0(\F)$-cash invariant} if $\rho(x+\eta)=\rho(x)-\eta$ for all $\eta\in L^0(\F)$.
\end{enumerate}
We say that $\rho$ is a \emph{conditional risk measure} on $L^\phi_{\F}(\mathcal{E})$ if $\rho$ is monotone, $L^0(\F)$-convex and $L^0(\F)$-cash invariant.  
The \emph{conjugate} of a conditional risk measure $\rho$ on $L^\phi_\F(\mathcal{E})$ is defined as
\begin{equation}
\label{eq:FenchelConjugate2}
{\rho}^\#(y):=\esssup\{\E[x y|\F]-\rho(x) \colon x\in L^\phi_\F(\mathcal{E})\} \quad \mbox{ for }y\in L^{\phi^\ast}_\F(\mathcal{E}).
\end{equation}
\end{definition}

\begin{lemma}\label{lem:localProp}
Let $\rho$ be a conditional risk measure on $L^\phi(\mathcal{E})$.  
Then $1_A\rho(x)=1_A\rho(1_A x)$ for all $x\in L^\phi(\EE)$ and $A\in\F$. 
\end{lemma}
\begin{proof}

Given $A\in\F$ and $x\in L^\phi(\mathcal{E})$, by $L^\infty(\F)$-convexity, one has 
\begin{align*}
\rho(1_A x)&=\rho(1_A x + 1_{A^c} 0)\\
&\le 1_A\rho(x)+1_{A^c}\rho(0)\\
&=1_A\rho(1_A (1_A x) + 1_{A^c}x)+1_{A^c}\rho(0)\\
&\le 1_A\rho(1_A x)+1_{A^c}\rho(0).
\end{align*} 
Multiplying by $1_A$, all the inequalities become equalities and the conclusion follows.
\end{proof}

The proof of the following result is easy to find and therefore omitted.
\begin{lemma}\label{lem:supStHull}
Let $\emptyset\neq S\subset L^0(\mathcal{E})$. 
Then $\esssup S=\esssup\left( s_\F(S)\right)$.
\end{lemma}

\begin{proposition}\label{prop:extension}
Every conditional risk measure $\rho$ on $L^\phi(\mathcal{E})$   extends to a unique conditional risk measure $\bar{\rho}$ on $L^\phi_\F(\mathcal{E})$.  Moreover, the extension $\bar{\rho}$  satisfies  $
\bar{\rho}^\#(y)=\rho^\#(y)$ for all $y\in L^{\phi^\ast}(\mathcal{E})$.
\end{proposition}
\begin{proof}
Given $x\in L^{\phi}_\F(\mathcal{E})$, due to Proposition~\ref{prop: ccHull}, there exist a sequence $(x_k)_{k\in\N}\subset L^\phi(\mathcal{E})$ and a countable partition $(A_k)_{k\in\N}\subset\F$ of $\Omega$ such that $x=\sum_{k\in\N} 1_{A_k}x_k$. 
Set $\bar{\rho}(x):=\sum_{k\in\N} 1_{A_k}\rho(x_k)$. 
Lemma \ref{lem:localProp} allows to show that the mapping $\bar{\rho}\colon L^\phi_\F(\mathcal{E})\to L^0(\F)$ is well-defined.  
Let us prove that $\bar{\rho}$ is $L^0(\F)$-cash invariant. 
Fix $x\in L^{\phi}_\F(\mathcal{E})$ and $\eta\in L^0(\F)$. 
Put $x=\sum_{k\in\N} 1_{A_k}x_k$ where $(x_k)\subset L^\phi(\mathcal{E})$ and $(A_k)\subset\F$ is a countable partition of $\Omega$.  
For $k=1,2,\cdots$ define $B_k:=\{k-1\le|\eta|<k\}$. 
Set $C_{i,j}:=A_i\cap B_j$ for all $i,j\in\mathbb{N}$. 
Then, $1_{C_{i,j}}\eta\in L^\infty(\F)$ and $1_{C_{i,j}}x\in L^\phi(\mathcal{E})$ for all $i,j$. 
Using Lemma \ref{lem:localProp} and $L^\infty(\F)$-cash invariance we have
\begin{align*}
\bar{\rho}(x+\eta)&=\sum_{i,j\in\N} 1_{C_{i,j}}\rho\big(1_{C_{i,j}}(x+\eta)\Big)\\
&=\sum_{i,j\in\N} 1_{C_{i,j}}\rho\big(1_{C_{i,j}}x+1_{C_{i,j}}\eta\Big)\\
&=\sum_{i,j\in\N} 1_{C_{i,j}}\Big(\rho(1_{C_{i,j}} x)- 1_{C_{i,j}}\eta\Big)\\
&=\sum_{i,j\in\N} 1_{C_{i,j}}\rho(1_{C_{i,j}} x)- \sum  1_{C_{i,j}}\eta\\
&=\bar{\rho}(x)- \eta.
\end{align*}
A similar argument shows that $\bar{\rho}$ is $L^0(\F)$-convex. 
The extension $\bar{\rho}$ is unique. 
Indeed, suppose that $\tilde{\rho}$ is another conditional risk measure on $L^\phi_\F(\EE)$ extending  $\rho$. 
Fix $x\in L^\phi_\F(\EE)$, and write $x=\sum_{k\in\N} 1_{A_k}x_k$ where $(x_k)_{k\in\N}\subset L^\phi(\mathcal{E})$ and $(A_k)_{k\in\N}\subset\F$ is a countable partition of $\Omega$. 
Due to \cite[Theorem 3.2]{kupper03} both $\bar{\rho}$ and $\tilde{\rho}$ have the local property.  
Then, since $\bar{\rho}$ and $\tilde{\rho}$ agree on $L^\phi(\EE)$, we get
\[
\bar{\rho}(x)=\sum_{k\in\N} 1_{A_k} \bar{\rho}(1_{A_k}x_k)=\sum_{k\in\N} 1_{A_k} \bar{\rho}(x_k)=\sum_{k\in\N} 1_{A_k} \tilde{\rho}(x_k)=\sum_{k\in\N} 1_{A_k} \tilde{\rho}(1_{A_k}x_k)=\tilde{\rho}(x). 
\]  
Finally, Lemma \ref{lem:supStHull} yields that $\rho^\#(y)=\bar{\rho}^\#(y)$ for all $y\in L^{\phi^\ast}(\mathcal{E})$. 
\end{proof}


\begin{definition}
Let $f\colon L^{\phi^\ast}_\F(\mathcal{E})\to\bar{L}^0(\F)$ be a function.  
For every $\eta\in L^0(\F)$, we define the sublevel set 
$$V_\eta(f):=\left\{ y\in L^{\phi^\ast}_\F(\mathcal{E}) \colon f(y)\le \eta \right\}.$$
We say that $f$ is stably  $w^\#_s$-inf-compact if  $V_\eta(f)$ is either stably $w^\#_s$-compact or empty for every $\eta\in L^0(\F)$.
\end{definition}

Next, we state the main result of this section. 
\begin{theorem}\label{th:dualRep}
Let $\rho$ be a conditional risk measure on $L^\phi(\mathcal{E})$ and with values  in $L^1(\F)$,  which has the representation
\begin{equation}\label{eq:rep}
\rho(x)=\underset{y\in L^{\phi^\ast}(\EE)}\esssup\left\{ \mathbb{E}[x y|\F]-\rho^\#(y)\right\}
\quad\mbox{ for all }x\in L^\phi(\mathcal{E}).
\end{equation}
Then, the following conditions are equivalent:
\begin{enumerate}
\item $\rho$ has the Lebesgue property: i.e. 
		\[
	y\in L^\phi(\mathcal{E}),\:|x_n|\leq |y|\mbox{ for all }n\in\N,\:x_n\to x\text{ a.s. implies }\rho(x_n)\to\rho(x)\mbox{ a.s.}.
	\]
	\item $\bar{\rho}$ has the Lebesgue property: i.e. 
		\[
	y\in L^\phi_{\F}(\mathcal{E}),\:|x_n|\leq |y|\mbox{ for all }n\in\N,\:x_n\to x\mbox{ a.s. implies }\bar{\rho}(x_n)\to\bar{\rho}(x)\mbox{ a.s.}.
	\]
	\item For each $x\in L^\phi(\mathcal{E})$ there exists $y\in L^{\phi^\ast}_\F(\mathcal{E})$ with $y\le 0$ and $\mathbb{E}[y|\F]=-1$ such that $\rho(x)=\mathbb{E}[ xy|\F]-\bar{\rho}^\#(y)$.
\item For each $x\in L^\phi_\F(\mathcal{E})$ there exists $y\in L^{\phi^\ast}_\F(\mathcal{E})$ with $y\le 0$ and $\mathbb{E}[y|\F]=-1$ such that $\bar{\rho}(x)=\mathbb{E}[ xy|\F]-\bar{\rho}^\#(y)$.
	
	\item $\bar{\rho}^\#$ is stably $w^\#_s$-inf-compact.
\end{enumerate}
\end{theorem}

\begin{remark}\label{rem:rep}
If $\rho$ admits the representation \eqref{eq:rep}, then by Proposition~\ref{prop:extension} it follows that
\begin{equation}\label{eq:repII}
\bar{\rho}(x)=\underset{y\in L^{\phi^\ast}_{\F}(\EE)}\esssup\left\{ \mathbb{E}[x y|\F]-\bar{\rho}^\#(y)\right\}
\quad\mbox{ for all }x\in L^\phi_{\F}(\mathcal{E}).
\end{equation}
In fact, both \eqref{eq:rep} and \eqref{eq:repII} are equivalent. 
In the forthcoming discussion, we also  show that \eqref{eq:repII} is also equivalent to the lower semicontinuity of $\rho$; see Remark~\ref{rem:lsc} below. 
\end{remark}

Before proving Theorem~\ref{th:dualRep}, we need to introduce some terminology and give some preliminary results.
\begin{definition}\label{defn:lsc}
Let $f\colon L^{\phi}_\F(\mathcal{E})\to\bar{L}^0(\F)$ be a function.
\begin{itemize}
\item We say that $f$ is $w_s$-lower semicontinuous if for every $\varepsilon\in L^0_{++}(\F)$ and $x\in L^\phi_\F(\mathcal{E})$ there exists a $w_s$-neighborhood $U$ of $x$ such that 
 $f(y)\ge f(x)-\varepsilon$ for all $y\in U$.  
\item We say that $f$ is $w_s$-upper semicontinuous if $-f$ is  $w_s$-lower semicontinuous. 
\end{itemize}
We define analogous notions for the $w^\#_s$-topology. 
\end{definition}

\begin{remark}
The definition of lower (and upper) semi continuous function given in the present paper differs from that given in \cite{kupper03,jamzap}. 
Namely, in \cite{kupper03,jamzap}, the function $f\colon L^{\phi^\ast}_{\F}(\EE)\to \bar{L}^0(\F)$ is called lower semicontinuous if the sublevel set $V_\eta(f)$ is closed for every $\eta\in\bar{L}^0(\F)$. 
We claim that if $f$ is lower semicontinuous in the sense of Definition~\ref{defn:lsc}, then $f$ is also lower semicontinuous in the sense of \cite{jamzap}. 
To prove that, suppose that $f$ is lower semicontinuous as defined in Definition~\ref{defn:lsc} and fix $\eta\in\bar{L}^0(\F)$. 
If $x\notin V_\eta(f)$, then there exists $A\in\F$ with $\P(A)>0$ such that $f(x)>\eta$ on $A$. 
Let $\varepsilon\in L^0_{++}(\F)$ satisfy  $f(x)>f(x)-\varepsilon\ge\eta$ on $A$. 
Then, there exists a neighborhood $U$ of $x$ such that $f(y)\ge f(x)-\varepsilon$ for all $y\in U$. 
In particular, $f(y)>\eta$ on $A$ for all $y\in U$, which shows that $U\subset V_\eta(f)^c$. 
Since $x$ was arbitrary, we conclude that $V_\eta(f)$ is closed. 
\end{remark}

\begin{proposition}\label{prop:supf}
Let $(f_i)_{i\in I}$ be a family of functions from $L^\phi_\F(\mathcal{E})$ to $\bar{L}^0(\F)$.   
If $f_i$ is $w_s$-lower semicontinuous and has the local property for all $i\in I$, then the function $x\mapsto\esssup_{i\in I} f_i(x)$ is $w_s$-lower semicontinuous.
\end{proposition}
\begin{proof}
Fix $x\in L^{\phi}_\F(\mathcal{E})$ and $\varepsilon\in L^0_{++}(\F)$. 
For every $i\in I$, there exists an open neighborhood $U_i$ of $x$ such that $f_i(y)\ge f_i(x)-\varepsilon$ for all $y\in U_i$. 
Due to Lemma \cite[Lemma 2.4]{zapata}, there exists a countable partition $(A_k)_{k\in\N}\subset\F$ of $\Omega$ and a sequence $(i_k)_{k\in\N}\subset I$ such that $\sum_{k\in\N}1_{A_k} f_{i_k}(x)\ge \esssup_{i\in I} f(x) -\varepsilon$. 
Take $U:=\sum_{k\in \N} 1_{A_{k}} U_{i_k}$, which is an open neighborhood of $x$ (see Remark \ref{rem:stTop}). 
Then, since $f_{i_k}$ has the local property it follows that 
\[
\esssup_{i\in I} f_i(y)\ge  \sum_{k\in\N}1_{A_k} f_{i_k}(y)\ge \sum_{k\in\N}1_{A_k} f_{i_k}(x)-\varepsilon  \ge \esssup_{i\in I} f_i(x) -\varepsilon 
\]
for every $y\in U$. 
The proof is complete. 
\end{proof}

\begin{remark}\label{rem:lsc}
 The function $x\mapsto \mathbb{E}[x y|\F]$ is $w_s$-continuous.    
 As a consequence of Proposition \ref{prop:supf}, if $\bar{\rho}\colon L^\phi_\F(\EE)\to L^0(\F)$ has the representation~\eqref{eq:repII} (or equivalently \eqref{eq:rep} holds), then $\bar{\rho}$ is $w_s$-lower semicontinuous.   
 Conversely, the representation \eqref{eq:repII} holds true whenever $\rho$ is $w_s$-lower semicontinuous, which follows from the extended Fenchel-Moreau theorem~\cite[Theorem 3.8]{kupper03}.  
\end{remark}

\begin{lemma}\label{lem:dom}
Let $\rho$ be a conditional risk measure on  $L^{\phi}_{\F}(\EE)$. 
If ${\rho}^\#(y)<+\infty$ for $y\in L^{\phi^\ast}_{\F}(\EE)$, then $y\le 0$ and $\mathbb{E}[y|\F]=-1$. 
\end{lemma}
\begin{proof}
Fix $\lambda\in L^0_{++}(\F)$. 
By monotonicity, we have
\[
+\infty>{\rho}^\#(y)\ge \E[\lambda 1_{\{y\geq 0\}} y|\F]-{\rho}(\lambda 1_{\{y\geq 0\}})\ge \lambda \E[y^+|\F] -{\rho}(0).
\]
Since $\lambda$ is arbitrary, we have that the left side of the inequality above is finite only if $y\le 0$.
Moreover, given $\lambda\in L^0_{++}(\F)$, by $L^0(\F)$-cash invariance, we have 
\[
+\infty>{\rho}^\#(y)\geq \E[\lambda y|\F]-{\rho}(\lambda)=\lambda \left(\E[y|\F]+1\right) -{\rho}(0).
\]
Since $\lambda$ is arbitrary, necessarily $\E[y|\F]=-1$.
\end{proof}

A similar result can be found in \cite[Lemma 2]{kabanov}.
\begin{lemma}\label{lem:kabanov}
\label{lem: seq}
Let $(y_n)_{n\in\N}$ be a sequence in $L^0(\F)$ such that $\limsup_{n\to+\infty} y_n=y$. 
Then, there exists a sequence $\mathfrak{n}_1<\mathfrak{n}_2<\cdots$ in $L^0(\F,\N)$ such that $y_{\mathfrak{n}_k}\to y$~a.s.. 
\end{lemma}

We need the following scalarization lemma. 
\begin{lemma}\label{lem:scalarization}
Let $\rho$ be a conditional risk measure on $L^\phi(\EE)$ and with values in $L^1(\F)$.  
Consider the function $\rho_0\colon L^\phi(\EE)\to \R$ defined by $\rho_0(x):=\E[\rho(x)]$, and its convex conjugate $\rho_0^{\#}\colon L^{\phi}(\EE)\to \R\cup\{+\infty\}$ defined by  $\rho_0^\#(y)=\sup_{x\in L^{\phi}(\EE)}\{\E[x y]-\rho_0(x)\}$.  
Then, 
\begin{enumerate}
\item $\rho^\#_0(y)=\mathbb{E}[\rho^\#(y)]$ for all $y\in L^{\phi^\ast}(\EE)$.
\item If $\rho$ has the representation \eqref{eq:rep}, 
then 
\begin{equation}
\label{eq:reprho0}
\rho_0(x)=\underset{y\in L^{\phi^\ast}(\EE)}\sup\left\{ \mathbb{E}[x y]-\rho^\#_0(y)\right\}
\quad 
\mbox{ for all }x\in L^\phi(\EE). 
\end{equation}
In that case, $\rho_0$ is $\sigma(L^{\phi}(\EE),L^{\phi^\ast}(\EE))$-lower semicontinuous.
\end{enumerate}
\end{lemma}
\begin{proof}
The proof is a closed to the  arguments used in the proof of~\cite[Theorem 1]{DetlefsenScandolo}. 
 Fix $y\in L^{\phi^\ast}(\EE)$. 
 The set $M_y:=\left\{\mathbb{E}[x y|\F]-\rho(x) \colon x\in L^{\phi}(\EE)\right\}$ is upward directed. 
 Therefore, there exists a sequence $(x_n)_{n\in\N}\subset L^{\phi}(\EE)$ such that 
 $\mathbb{E}[x_n y|\F]-\rho(x_n)\uparrow \esssup M_y=\rho_0^\#(y)$ a.s.. 
 Notice that, since $\rho^\#(y)\ge -\rho(0)\in L^1(\F)$, the expectation $\mathbb{E}[\rho^\#(y)]$ exists (allowing the value $+\infty$). 
 By monotone convergence, we have
 \[
 \mathbb{E}[\rho^\#(y)]=\mathbb{E}\Big[\underset{n\to+\infty}\lim(\mathbb{E}[x_n y|\F]-\rho(x)) \Big]
 = \underset{n\to+\infty}\lim (\mathbb{E}[x_n y]-\rho_0(x))\le \rho^\#(y).
 \]
 If $x\in L^\phi(\EE)$, we have that $\rho^\#(x)\ge \mathbb{E}[xy|\F]-\rho(x)$, hence
 \[
 \mathbb{E}[\rho^\#(y)]\ge\mathbb{E}\Big[  \mathbb{E}[xy|\F]-\rho(x) \Big]
 =\mathbb{E}[xy]-\rho_0(x). 
 \]
 Then, taking the essential supremum over all $x\in L^\phi(\EE)$, we get $\mathbb{E}[\rho^\#(y)]\ge \rho^\#_0(y)$. 
 We conclude that $\rho^\#_0(y)=\mathbb{E}[\rho^\#(y)]$. 
 Finally, if $\rho$ has the representation \eqref{eq:rep}, a similar argumentation shows that $\rho_0$ is represented as in~\eqref{eq:reprho0}. 
 In that case, $\rho_0$ is the supremum of lower semicontinuous functions, hence it is also lower semicontinuous.  
\end{proof}

The following result is an adaptation to the present setting of \cite[Lemma 2.3]{owari2014lebesgue}. 
\begin{lemma}
\label{lem: ineq}
Let $\rho$ be a conditional risk measure on $L^{\phi}(\mathcal{E})$ with $\rho(0)=0$ which admits the representation \eqref{eq:rep}.  
Then, for any $\beta\in L^0(\F)$, $x\in L^{\phi}_\F(\mathcal{E})$ and $y\in L^{\phi^\ast}_\F(\mathcal{E})$ with $y\le 0$, 
\[
\E[x y|\F]-\bar{\rho}^\#(y)\ge -\beta \quad\mbox{implies}\quad \bar{\rho}^\#(y)\le 2\beta + 2\bar{\rho}(-2|x|).
\]
\end{lemma}
\begin{proof}
Due the representation \eqref{eq:repII}, we have
\[
0=\bar{\rho}(0)=-\underset{y\in L^{{\phi}^\ast}_\F(\mathcal{E})\colon y\le 0}\esssup \bar{\rho}^\#(y),
\]
where the essential supremum can be taken over $y\le 0$ due to Lemma \ref{lem:dom}. 
Then, given $n\in\N$, by \cite[Lemma 2.4]{zapata} there exists $y_n\in L^{{\phi}^\ast}_\F(\mathcal{E})$ with $y_n\le 0$ such that $0\leq\bar{\rho}^\#(y_n)\le 1/n$. 
Thus, 
\begin{align*}
\E[x y|\F]\le&\E\left[-2|x|\frac{y+y_n}{2}|\F\right]\leq \bar{\rho}(-2|x|) + \bar{\rho}^\ast\left(\frac{y+y_n}{2}\right)\\
\le& \bar{\rho}(-2|x|) + \frac{1}{2}\bar{\rho}^\ast(y) + \frac{1}{2}\bar\rho^\ast(y_n)=\bar{\rho}(-2|x|) + \frac{1}{2}\bar{\rho}^\ast(y) + \frac{1}{n}.
\end{align*}
Letting $n\to+\infty$, we obtain
\[
\E[x y|\F]\leq\bar{\rho}(-2|x|) + \frac{1}{2}\bar{\rho}^\#(y).
\]
Then, if $\E[x y|\F]-\bar{\rho}^\#(y)\ge -\beta$, we have that
\[
\bar{\rho}^\#(y)\leq \E[x y|\F]+\beta\leq \bar{\rho}(-2|x|) + \frac{1}{2}\bar{\rho}^\#(y) +\beta. 
\] 
We finally get $
\bar{\rho}^\#(y)\leq 2\beta + 2\bar{\rho}(-2|x|)$ and the proof is complete.
\end{proof}

We state \cite[Theorem 5.13]{jamzap} in the present setting. 
\begin{lemma}\label{lem: upper}
Let $S\subset L^{\phi^\ast}_{\F}(\EE)$ be stably $w^\#_s$-compact. 
Suppose that $f\colon L^{\phi^\ast}_{\F}(\EE)\to \bar{L}^0(\F)$ is $w^\#_s$-upper semicontinous, has the local property, and $f(S)\subset L^0(\F)$. 
Then, there exists $y_0\in S$ such that $f(y_0)=\esssup_{y\in S} f(y)$.
\end{lemma}

We now turn to the proof of Theorem~\ref{th:dualRep}.

\begin{proof}
$(ii)\Rightarrow(i)$ and $(iv)\Rightarrow(iii)$ are clear.

$(i)\Rightarrow (iii):$ Let $\rho_0\colon L^\phi(\EE)\to \R$ be defined by $\rho_0(x):=\E[\rho(x)]$.  
Fix $x\in L^\phi(\EE)$. 
Due to \eqref{eq:rep} we have that $\rho(x)\ge \E[xy|\F]-\rho^\#(y)$ for all $y\in L^{\phi^\ast}(\EE)$. 
It suffices to show that there exists $y\in L^{\phi^\ast}(\EE)$ such that 
\[
\mathbb{E}[\rho(x)]=\E\left[\E[xy|\F]-\rho^\#(y)\right],
\]
or equivalently that $\rho_0(x)=\mathbb{E}[x y]-\rho^\#_0(y)$ in view of Lemma~\ref{lem:scalarization}. 
Due to the representation \eqref{eq:rep} and Lemma~\ref{lem:scalarization}, we have that $\rho_0$ is $\sigma(L^\phi(\EE),L^{\phi^\ast}(\EE))$-lower semicontinuous. 
 By \cite[Theorem 1.1]{owari2014lebesgue} it suffices to prove that $\rho_0$ has the Lebesgue property. 
 Indeed, suppose that $x_n\to x$ a.s. and $|x_n|\le y$ for some $y\in L^{\phi^\ast}(\EE)$. 
 Since $\rho$ has the Lebesgue property, it holds that $\rho(x_n)\to\rho(x)$ a.s.. 
 By monotonicity, we have that $|\rho(x_n)|\le |\rho(y)|\vee |\rho(-y)|$. 
 Then, by dominated convergence we obtain
 \[
 \underset{n\to+\infty}\lim \rho_0(x_n)
 =\E\left[\underset{n\to+\infty}\lim \rho(x_n) \right]
 =\E[\rho(x)]
 =\rho_0(x).
 \]   

$(iii)\Rightarrow (i):$ Suppose that $x_n\to x$ a.s. and $|x_n|\le y$ for some $y\in L^{\phi}(\EE)$. 
Then, by dominated convergence we have that $\lim_{n\to+\infty}\E[x_n z|\F]=\E[x z|\F]$ for all $z\in L^{\phi^\ast}(\mathcal{E})$. 
Then,
\begin{align*}
\rho(x)=&\underset{z\in L^{{\phi}^\ast}(\mathcal{E})}\esssup\left\{\underset{n\to+\infty}\lim\E[x_n z|\F]-{\rho^\#}(z)\right\}\\
\le & \underset{n\to+\infty}\liminf\,\underset{z\in L^{{\phi}^\ast}(\mathcal{E})} \esssup\left\{ E[x_n z|\F]-{\rho^\#}(z)\right\}\\
=&\underset{n\to+\infty}\liminf\rho(x_n).
\end{align*}
It suffices to show that $\rho(x)\ge \limsup_{n\to+\infty}\rho(x_n)$. 
By contradiction, suppose that $\rho(x)<\limsup_{n\to+\infty}\rho(x_n)$ on some $A\in\F$ with $\P(A)>0$. 
Since $\bar{\rho}$ has the local property, we can  suppose w.l.o.g. that $A=\Omega$ by considering the restrictions of random variables and traces of $\sigma$-algebras on $A$.    
We consider $\rho_0(x):=\E[\rho(x)]$, for which we have 
\[
\rho_0(x)< \E\left[\limsup_{n\to+\infty}\rho(x_n)\right].
\] 
By Lemma~\ref{lem:kabanov}, there exists $\mathfrak{n}_1<\mathfrak{n}_2<\cdots$ in $L^0(\F;\N)$ such that  
 $\lim_{k\to+\infty}\rho(x_{\mathfrak{n}_k})=\limsup_{n\to+\infty}\rho(x_n)$ a.s.. 
 For each $k$, take $z_k:=x_{\mathfrak{n}_k}$. 
 We have that $|z_k|\leq|y|$ for all $k\in\N$, and $\lim_{k\to+\infty} z_k=x$ a.s..  
Then, by dominated convergence
\begin{equation}\label{eq:LepProbIne}
\lim_{k\to+\infty}\rho_0(z_k)=\E\left[\lim_{k\to+\infty}\rho(z_k)\right]= \E\left[\limsup_{n\to+\infty}\rho(x_n)\right]>\rho_0(x).
\end{equation}
On the other hand, by assumption, for each $z\in L^{\phi}(\mathcal{E})$ there exists  $y\in L^{{\phi}^\ast}(\mathcal{E})$ such that
\[
\rho_0(z)=\E\left[\E[z y|\F] - \rho^\#(z)\right]=\E[z y] - \rho^\#_0(y).
\] 
Thus, by \cite[Theorem 1.1]{owari2014lebesgue}, we have that $\rho_0$  necessarily has the Lebesgue property; hence, $\lim_{k\to+\infty}\rho_0(z_k)=\rho_0(x)$, which is a contradiction to \eqref{eq:LepProbIne}.

$(iii)\Rightarrow(iv):$ If $x\in L^{{\phi}}_\F(\mathcal{E})$, then $x=\sum_{k\in\N} 1_{A_k}x_k$ for some sequence $(x_k)_{k\in\N}\subset L^{{\phi}}(\mathcal{E})$ and countable partition $(A_k)_{k\in\N}\subset\F$ of $\Omega$. 
For every $k\in\N$, take $y_k\in L^{{\phi}^\ast}(\mathcal{E})$ such that $\rho(x_k)=\mathbb{E}[x_k y_k]-\rho^\#(y_k)$. 
Put $y=\sum_{k\in\N} 1_{A_k}y_k$. 
Then, by lemma~\ref{lem:localProp},
\[
\rho(x)=\sum_{k\in\N} 1_{A_k}(\mathbb{E}[x_k y_k]-\rho^\#(y_k))
=\sum_{k\in\N} 1_{A_k}(\mathbb{E}[x_k y_k]-\bar{\rho}^\#(y_k))
=\mathbb{E}[x y]-\bar{\rho}^\#(y),
\]
where we have used that $\bar{\rho}^\#(y_k)={\rho}^\#(y_k)$. 

$(i)\Rightarrow(ii)$: Suppose that $x_n\rightarrow x$ a.s. in $L^\phi_\F(\EE)$ with $|x_n|\leq y$ for some $y\in L^{{\phi}}_\F(\mathcal{E})$. 
Due to Proposition \ref{prop: ccHull}, we can take $(y_k)_{k\in\N}\subset L^{{\phi}}(\mathcal{E})$ and a countable partition $(A_k)_{k\in\N}\subset\F$ of $\Omega$ such that $y=\sum_{k\in\N} 1_{A_k} y_k$. 
Then, for each $k\in\N$, $|1_{A_k} x_n|\leq 1_{A_k} y_k$ and $1_{A_k}x_n\rightarrow 1_{A_k} x$ a.s.. 
Thus, $\rho(1_{A_k} x_n)\to \rho(1_{A_k} x)$ a.s. and in turn  $\bar{\rho}(x_n)=\sum_{k\in\N} 1_{A_k} \rho(1_{A_k}x_n)$ converges a.s. to $\sum 1_{A_k} \rho(1_{A_k}x)=\bar{\rho}(x)$.

$(i)\Rightarrow(v)$: 
We have already seen in $(i)\Rightarrow(iii)$ that $\rho_0$ has the Lebesgue property. 
Then, due to \cite[Theorem 1.1]{owari2014lebesgue}, we have that, for any $c\in\R$, the set 
\[
V^0_c:=\left\{y\in L^{{\phi}^\ast}(\mathcal{E})\colon \rho^\#_0(y)\leq c \right\} \quad\mbox{ is }\sigma(L^{{\phi}^\ast}(\mathcal{E}),L^{{\phi}}(\mathcal{E}))\mbox{-compact.}
\]
We next prove that 
\[
V_c:=V_c(\bar{\rho}^\#):=\left\{y\in L^{{\phi}^\ast}_{\F}(\mathcal{E})\colon \bar{\rho}^\#(y)\leq c \right\}
\]
  is stably $w^\#_s$-compact. 
Define
\[
f_c(x):=\esssup\left\{\E[|x{y}||\F]\colon y\in {L}^{{\phi}^\ast}_\F(\mathcal{E}),\:\bar{\rho}^\#(y)\leq c\right\} \quad \mbox{ for }{x}\in{L}^{{\phi}}_\F(\mathcal{E}),
\]
\[
f_{0,c}(x):=\sup\left\{\E[|x y|]\colon y\in {L}^{{\phi}^\ast}(\mathcal{E}),\:\rho_0^\#(y)\leq c\right\} \quad \mbox{ for }x\in L^{\phi}(\mathcal{E}).
\]
Fix $x\in L^{\phi}(\mathcal{E})$. 
Since the functions $y\mapsto \E[|x{y}||\F]$ and $\bar{\rho}$ have the local property, the set $H_{c}(x):=\left\{\E[|x{y}||\F]\colon y\in{L}^{{\phi}^\ast}_{\F}(\mathcal{E}),\:\bar{\rho}^\#(y)\leq c\right\}$ is upward directed.  
We can find a sequence $(y_n)\subset {L}^{{\phi}^\ast}_{\F}(\mathcal{E})$ with ${\bar{\rho}}^\#(y_n)\leq c$ such that $\E[|x{y_n}||\F]\uparrow \esssup H_c(x)$ a.s.. 
 Since $c\geq\E[{\bar{\rho}}^\#(y_n)]={\rho}_0^\#(y_n)$ for each $n\in\N$, by monotone convergence, we have that $\E[f_c(x)]=\E[\lim_{n\to+\infty} \E[|x{y_n}||\F]]=\lim_{n\to+\infty} \E[|x{y_n}|]\leq f_{0,c}(x)$. 
 We conclude that
\begin{equation}
\label{eq: fc}
\E[f_c(x)]\leq f_{0,c}(x)\quad\textnormal{ for all }x\in L^{\phi}(\mathcal{E}). 
\end{equation}
Now, suppose that $x_n\downarrow 0$ a.s. in $L^{\phi}_\F(\mathcal{E})$. 
Due to Proposition~\ref{prop: ccHull}, there is a countable partition $(A_k)_{k\in\N}\subset\F$ of $\Omega$ such that $1_{A_k}x_1\in L^{\phi}(\mathcal{E})$ for all $k\in\N$. 
Since $L^{\phi}(\mathcal{E})$ is solid and $0\leq x_n\leq x_1$, we have that $1_{A_k}x_n\in L^{\phi}(\mathcal{E})$ for all $k,n\in\N$.    

Since $V_c^0$ is $\sigma({L}^{{\phi}^\ast}(\mathcal{E}),{L}^{{\phi}}(\mathcal{E}))$-compact, it follows by \cite[Lemma 2.1]{owari2014lebesgue} that\newline $\lim_{n\to+\infty}f_{0,c}(1_{A_k}x_n)=0$ a.s. for each $k\in\N$. 
Therefore, in view of \eqref{eq: fc}, one has 
\[
0\leq\lim_{n\to+\infty}\E[f_c(1_{A_k}x_n)]\leq \lim_{n\to+\infty} f_{0,c}(1_{A_k}x_n)=0\mbox{ a.s. for all }k\in\N. 
\]

Since $f_c(1_{A_k}x_1)\ge f_c(1_{A_k}x_{2}) \ge\cdots\ge0$, necessarily $f_c(1_{A_k}x_n)\rightarrow 0$ a.s. for each $k\in\N$. 
Since $f_c$ has the local property, we conclude that $f_c(x_n)\rightarrow 0$ a.s.
By Theorem~\ref{thm:criterion}, we obtain that ${V}_c$ is stably $w_s^\#$-compact (notice that $V_c$ is $w^\#_s$-closed as $\bar{\rho}^\#$ is $w_s^\#$-lower semicontinuous).   

Finally, if $\eta\in L^0(\F)$ is arbitrary, we can find a sequence $(c_k)_{k\in\N}\subset\R$ and a countable partition $(A_k)_{k\in\N}\subset \F$ of $\Omega$ such that $\sum_{k\in\N}1_{A_k}c_k\le \eta$. 
Then, $V_\eta(\rho^\#)\subset \sum_{k\in\N} 1_{A_k}V_{c_k}$. 
We have that $\sum_{k\in\N} 1_{A_k}V_{c_k}$ is stably $w_s^\#$-compact, hence $V_\eta(\rho^\#)$ enjoys the same property. 

$(v)\Rightarrow(iv)$: 
We can suppose w.l.o.g. that $\rho(0)=0$. 
Due to the $w_s^\#$-lower semicontinuity of $\bar{\rho}^\#$, the function $y\mapsto\E[{x}{y}|\F]-\bar{\rho}^\#(y)$ is $w_s^\#$-upper semicontinuous.  
Then, the set 
\[
M_x:=\left\{y\in L^{{\phi}^\ast}_\F(\EE)\colon \E[{x}{y}|\F]-\bar{\rho}^\#(y)\ge \bar{\rho}(x)-1\right\}
\]
is  $w_s^\#$-closed and, %
applying Lemma \ref{lem: ineq} for $\beta=1-\bar{\rho}(x)$, 
we obtain that ${M}_{x}$ is contained in 
$V_{2-2{\bar{\rho}}({x})+{2}{\bar{\rho}}(-{2}|{x}|)}(\bar{\rho}^\#)$. 
Therefore, ${M}_{x}$ is stably $w_s^\#$-compact. 
Finally, due to Lemma~\ref{lem: upper} we obtain the result. 
\end{proof}

\begin{remark}
The second author of the present paper applied a model-theoretic method to obtain dual representation results for conditional risk measures on modules; see~\cite{zapataboolean,zapatathesis}. 
The approach in~\cite{zapataboolean,zapatathesis} relies on Boolean valued models, a tool from model-theory and logic  that formalizes  the Paul Cohen's method of forcing; see e.g.~\cite{kusraev2012boolean} for more details.
In \cite{zapataboolean,zapatathesis} it is shown that the transfer principle of Boolean-valued models allows for the systematic transcription of available duality results of classical risk measures as analogues for conditional risk measures. 
In particular, \cite[Theorem 3.6]{zapataboolean} essentially covers Theorem~\ref{th:dualRep} above. 
We point out that historically Theorem~\ref{th:dualRep} above appeared before than \cite[Theorem 3.6]{zapataboolean} in the arXiv version of the present paper; see~\cite{orihuela2}. 
Also, we believe that the approach given here is of independent interest  as it relies on conventional  methods which give an alternative path for working  mathematicians that are not familiar with the model-theoretic tools used in~\cite{zapataboolean}.     
\end{remark}

\end{document}